\documentclass[11pt]{article}
\usepackage{amsmath,amsthm}

\def\refeq#1{\if\workingver y(\ref{#1})-[[#1]]\else(\ref{#1})\fi}
\def\refth#1{\if\workingver y\ref{#1}-[[#1]]\else\ref{#1}\fi}
\def\mylabel#1{\if\workingver y\label{#1}{\bf\ \ [[#1]]\ \ }
\else\label{#1}\fi}
\def\mybibitem#1{\if\workingver y\bibitem{#1}{\bf\ \ [[#1]]\ \ }
\else\bibitem{#1}\fi}

\def\institute#1{\gdef\@institute{#1}}

\def\institutename{\par
 \begingroup
 \parskip=\z@
 \parindent=\z@
 \setcounter{@inst}{1}%
 \def\and{\par\stepcounter{@inst}%
 \noindent$^{\the@inst}$\enspace\ignorespaces}%
 \setbox0=\vbox{\def\thanks##1{}\@institute}%
 \ifnum\c@@inst=1\relax
 \else
   \setcounter{footnote}{\c@@inst}%
   \setcounter{@inst}{1}%
   \noindent$^{\the@inst}$\enspace
 \fi
 \ignorespaces
 \@institute\par
 \endgroup}

\newtheorem{thm}{Theorem}
\newtheorem{lem}[thm]{Lemma}

\newtheorem{cor}[thm]{Corollary}

\newtheorem{rem}[thm]{Remark}

\newtheorem{example}[thm]{Example}

\def\begeq#1{\begin{equation}\mylabel{#1}}
\def\endeq{\end{equation}}

\def\begalg{\begin{alg}}
\def\endalg{\end{alg}}

\let\workingver=n

\newcommand{\beq}{\begin{equation}}

\begin{document}

\title{Spectral Properties of a Binomial Matrix}
\author{
\vspace{1cm}
Pantelimon St\u anic\u a\thanks{Also Associated to the Institute of
Mathematics of Romanian Academy, Bucharest, Romania}
\ {\small and} Rhodes Peele\\
Auburn University Montgomery, Department of Mathematics,\\
Montgomery, AL 36117, USA\\
e-mail: stanpan@strudel.aum.edu\
}
\date{\today}
\maketitle

\pagestyle{myheadings}

\baselineskip=2\baselineskip

\section{Introduction}

In \cite{PS}, Peele and St\u anic\u a studied $n\times n$
matrices with the $(i,j)$
 entry the binomial coefficient $\binom{i-1}{j-1}$
 (matrix $L_n$), respectively
$\binom{i-1}{n-j}$ (matrix $R_n$)
 and derived many interesting results
on powers of these matrices.
In \cite{Stanica2000}, the author found that the same is true for a
much larger class of what he called {\em netted matrices}, namely
matrices with entries satisfying a certain type of recurrence among the
entries of all $2\times 2$ cells.
In this paper we continue the work in \cite{PS, Stanica2000}.
We find the generating function for all entries of $R_n^e$ on
each row or column and
we find all eigenvalues (with their multiplicities)
of $R_n$ modulo $3$ and $5$. We also
give a precise relation
between $R_n$ and $L_n$ with the help of a permutation matrix
and we prove the conjecture stated in \cite{PS}.

%\newpage

\section{Some Generating Functions}

We denote by $a_{i,j}=\binom{i-1}{n-j}$. We observe that they satisfy
the recurrence
\begin{equation}
\label{recurrence}
a_{i,j-1}=a_{i-1,j-1}+a_{i-1,j},
\end{equation}
which can be extended for $i\geq 0,j\geq 0$,
using the boundary conditions $a_{1,n}= 1,\, a_{1,j}=0,\, j\neq n$. We
shall use the following consequences of the boundary conditions
and recurrence \refeq{recurrence}: $a_{i,j}=0$ for $i+j\leq n$,
and $a_{i,n+1}=0, 1\leq i\leq n$.
In \cite{PS} it was proved that the entries of $R_n^e$ satisfy the recurrence
\begin{equation}
\label{gen_fib}
F_{e-1} a_{i,j}^{(e)}=F_{e} a_{i-1,j}^{(e)}+
F_{e+1} a_{i-1,j-1}^{(e)}-F_{e} a_{i,j-1}^{(e)},
\end{equation}
where $F_e$ is the Fibonacci sequence.
We were unable to find closed forms for {\em all} entries of $R_n^e$,
however
we found generating functions for the entries in each row and column of
 $R_n^e$.
 We use the following lemma, which can be proven easily once
 one guesses the two formulas
 (see also \cite{PS,Stanica2000}).
 \begin{lem}
 \label{row1_col1}
The elements on the first row
and the first column of $R_n^e$ are given by
\begin{equation}
\label{first_row}
\begin{split}
a_{1,j}^{(e)} &=\binom{n-1}{j-1} F_{e-1}^{n-j} F_e^{j-1}\\
a_{i,1}^{(e)} &=F_{e-1}^{n-i} F_e^{i-1}.
\end{split}
\end{equation}
 \end{lem}
We extend the tableau $a_{i,j}^{(e)}$, for $i\geq 0, j\geq 0$, using
 the recurrence \refeq{gen_fib} and Lemma \refth{row1_col1}.
\begin{thm}
The generating function for the $i$-th row of $R_n^e$ is
\[
r_i^{(e)}(x)=\sum_{j\geq 1} a_{i,j}^{(e)}\, x^{j-1}=
(F_e+F_{e+1} x)^{i-1} (F_{e-1}+F_e x)^{n-i}
\]
and, if $e>1$, the generating function for the $j$-th column of $R_n^e$ is
\begin{eqnarray*}
c_j^{(e)}(x)&=&\left(
\frac{F_{e+1} x-F_{e}}{F_{e-1}-F_e x}\right)^{j-1}
\frac{F_{e-1}^n}{F_{e-1}-F_e x}\left[1+
\sum_{s=1}^{j-1}\binom{n}{s}
\left( \frac{F_e(F_{e-1}-F_e x)}{F_{e-1}(F_{e+1} x-F_{e}) }\right)^s\right].
\end{eqnarray*}
and $c_j^{(1)}(x)=(1-x)^{j-1-n} x^{n-j}$.
\end{thm}
\begin{proof}
Multiplying \refeq{gen_fib} by $x^{j-1}$ and summing for $j\geq 1$, we get
\begin{eqnarray*}
&&
F_{e-1}\sum_{j\geq 1} a_{i,j}^{(e)}\, x^{j-1}+
F_{e}\sum_{j\geq 1} a_{i,j-1}^{(e)}\, x^{j-1}\\
&&\qquad \qquad\qquad\qquad\qquad\qquad =
F_{e}\sum_{j\geq 1} a_{i-1,j}^{(e)}\, x^{j-1}+
F_{e+1}\sum_{j\geq 1} a_{i-1,j-1}^{(e)}\, x^{j-1}
\Longleftrightarrow\\
&&
F_{e-1} r_i^{(e)}(x)+F_e x r_i^{(e)}(x)+F_e a_{i,0}^{(e)}=
F_e r_{i-1}^{(e)}(x)+F_{e+1} x r_{i-1}^{(e)}(x)+F_{e+1} a_{i-1,0}^{(e)}
\Longleftrightarrow\\
&&
(F_{e-1}+F_e x) r_i^{(e)}(x)-(F_{e}+F_{e+1} x) r_{i-1}^{(e)}(x)=
F_{e+1} a_{i-1,0}^{(e)}-F_e a_{i,0}^{(e)}\\
&&=
F_{e-1} a_{i,1}^{(e)}-F_e a_{i-1,1}^{(e)}
\stackrel{\text Lemma\ \refth{row1_col1}}{=}
F_{e-1} F_{e-1}^{n-i} F_e^{i-1}-F_e F_{e-1}^{n-i+1} F_e^{i-2}=0.
\end{eqnarray*}
Therefore,
\[
(F_{e-1}+F_e x) r_i^{(e)}(x)=(F_{e}+F_{e+1} x) r_{i-1}^{(e)}(x),
\]
which implies
\[
r_{i}^{(e)}(x)=\left(\frac{F_e+F_{e+1} x}{F_{e-1}+F_e x}
\right)^{i-1} r_1^{(e)}(x).
\]
The generating function for the first row of $R_n^e$ is
\begin{eqnarray*}
r_1^{(e)}(x)
&=&\sum_{j\geq 1} \binom{n-1}{j-1} F_{e-1}^{n-j} F_e^{j-1} x^{j-1}\\
&=& \sum_{s\geq 0} \binom{n-1}{s} F_{e-1}^{n-s-1} F_e^{s} x^{s}\\
&=&
F_{e-1}^{n-1} \sum_{s\geq 0} \binom{n-1}{s}
\left(\frac{F_e x}{F_{e-1}}  \right)^s\\
&=&
F_{e-1}^{n-1} \left( 1+\frac{F_e x}{F_{e-1}} \right)^{n-1}=
(F_{e-1}+F_e x)^{n-1}.
\end{eqnarray*}
Thus,
\[
r_i^{(e)}(x)=\left( \frac{F_e+F_{e+1} x}{F_{e-1}+F_e x} \right)^{i-1}
 (F_{e-1}+F_e x)^{n-1}= (F_e+F_{e+1} x)^{i-1}
(F_{e-1}+F_e x)^{n-i}.
\]
and the first claim is proved.

It is trivial to find the generating function of the columns of $R_n$.
Now, for $e>1$, multiplying \refeq{recurrence}
by $x^{i-1}$ and summing for
$i\geq 1$, we get
\begin{eqnarray*}
&&
F_{e-1}\sum_{i\geq 1} a_{i,j}^{(e)}\, x^{i-1}+
F_{e}\sum_{i\geq 1} a_{i,j-1}^{(e)}\, x^{i-1}\\
&&\qquad \qquad\qquad =
F_{e}\sum_{i\geq 1} a_{i-1,j}^{(e)}\, x^{i-1}+
F_{e+1}\sum_{i\geq 1} a_{i-1,j-1}^{(e)}\, x^{i-1}
\Longleftrightarrow\\
&&
F_{e-1} c_j^{(e)}(x)+F_e c_{j-1}^{(e)}(x)=
F_e x c_j^{(e)}(x)\\
&&\qquad +F_{e+1} x c_{j-1}^{(e)}(x)
+F_e a_{0,j}^{(e)}+F_{e+1} a_{0,j-1}^{(e)}
\Longleftrightarrow\\
&&
(F_{e-1}-F_e x) c_j^{(e)}(x)+(F_{e}-F_{e+1} x) c_{j-1}^{(e)}(x)\\
&&=
F_e a_{0,j}^{(e)}+F_{e+1} a_{0,j-1}^{(e)}=
F_{e-1} a_{1,j}^{(e)}+F_e a_{1,j-1}^{(e)}
\\
%\end{eqnarray*}
%\begin{eqnarray*}
&&
\stackrel{\text Lemma\ \refth{row1_col1}}{=}
F_{e-1}\binom{n-1}{j-1} F_{e-1}^{n-j} F_e^{j-1}+
F_e \binom{n-1}{j-2} F_{e-1}^{n-j+1} F_e^{j-2}\\
&&=
\binom{n}{j-1} F_{e-1}^{n-j+1} F_e^{j-1}.
\end{eqnarray*}
If $e>1$,  the generating function for the first column of $R_n^e$ is
\begin{eqnarray*}
 c_1^{(e)}(x)
 &=& \sum_{i\geq 1} F_{e-1}^{n-i} F_e^{i-1} x^{i-1}
 = F_{e-1}^{n-1} \sum_{s\geq 0} \left( \frac{F_e x}{F_{e-1}} \right)^s\\
&=& F_{e-1}^{n-1} \frac{1}{1-\frac{F_e x}{F_{e-1}}}=
 \frac{F_{e-1}^n}{F_{e-1}-F_e x}.
\end{eqnarray*}
It is not difficult to obtain that a recurrence of the form
\[
\alpha c_j+\beta c_{j-1}=u_{j-1},
\]
has the solution
\begin{equation}
\label{gen_rec}
\begin{split}
c_j
&= \left(\frac{-\beta}{\alpha}\right)^{j-2} c_1+\frac{1}{\alpha}
\left(u_{j-1}- \frac{\beta}{\alpha} u_{j-2}+\cdots+
\left(\frac{-\beta}{\alpha}\right)^{j-2} u_1\right)\\
&= \left(\frac{-\beta}{\alpha}\right)^{j-2} c_1+\frac{1}{\alpha}
\sum_{s=1}^{j-1} u_s \left(\frac{-\beta}{\alpha}\right)^{j-s-1}.
\end{split}
\end{equation}
Using \refeq{gen_rec} in the recurrence for $c_j^{(e)}(x)$, we get
\begin{eqnarray*}
c_j^{(e)}(x)&=&\left(
\frac{F_{e+1} x-F_{e}}{F_{e-1}-F_e x}\right)^{j-1} c_1^{(e)}(x)\\
&&+
\frac{1}{F_{e-1}-F_e x}\sum_{s=1}^{j-1}\binom{n}{s}
\left(\frac{F_{e+1} x-F_{e}}{F_{e-1}-F_e x}\right)^{j-s-1}
F_{e-1}^{n-s} F_e^s\\
&=&\left(
\frac{F_{e+1} x-F_{e}}{F_{e-1}-F_e x}\right)^{j-1}
\frac{F_{e-1}^n}{F_{e-1}-F_e x}\left(1+
\sum_{s=1}^{j-1}\binom{n}{s}
\left( \frac{F_e(F_{e-1}-F_e x)}{F_{e-1}(F_{e+1} x-F_{e}) }\right)^s\right).
\end{eqnarray*}
\end{proof}

\section{Characteristic Polynomials Modulo $3$}

 In this section, we are interested in the eigenvalues of $R_n$ modulo $p$,
 and their multiplicities.
\begin{lem}
\label{trace}
We have
\[
trace(R_n)=F_{n}.
\]
\end{lem}
\begin{proof}
We need to prove
 \[
\sum_{i=1}^n \binom{i-1}{n-i}\stackrel{n-i=k}{=}
\sum_{k=0}^{n} \binom{n-k-1}{k}=F_{n},
 \]
 which is a well-known relation \cite{GKP}.
\end{proof}

 Our main result of this section is
\begin{thm}
The characteristic polynomial of $R_n$ modulo $3$, say $p_n(x)$, is
\begin{eqnarray*}
 p_{4k}(x)&=&(2+x+x^2)^k (2+2x+x^2)^k=(1+x^4)^k\\
 p_{4k+1}(x)&=&2(1+x)^{2[\frac{k+1}{2}]}(2+x)^{2[\frac{k}{2}]+1} (1+x^2)^k\\
&=& -(x^4-1)^k (x-1)\ \text{if $k$ even or}\
-(x^4-1)^{k}(x+1)\ \text{if $k$ odd}\\
 p_{4k+2}(x)&=&(2+x+x^2)^{2[\frac{k+1}{2}]} (2+2x+x^2)^{2[\frac{k}{2}]+1}%\\
  \end{eqnarray*}
 \begin{eqnarray*}
&=&
(x^4+1)^{k}(x^2-x-1)\ \text{if $k$ even or}\
(x^4+1)^{k}(x^2+x-1)\ \text{if $k$ odd}\\
 p_{4k+3}(x)&=& 2(1+x)^{2[\frac{k}{2}]+1}
 (2+x)^{2[\frac{k+1}{2}]}(1+x^2)^{k+1}\\
&=& -(x^4-1)^k (x+1)(x^2+1)\ \text{if $k$ even or}\\
&&\quad -(x^4-1)^{k}(x-1)(x^2+1)\ \text{if $k$ odd}.
\end{eqnarray*}
\end{thm}

\begin{proof}
%We use a method of \cite{PS} and \cite{Strauss}.
First, we prove that $R_{2k}^4\equiv -I_{2k}\pmod 3$ and
$R_{2k+1}^4\equiv I_{2k+1}\pmod 3$. In fact a more general result is
true; using Theorem 12 of \cite{PS}
we know that if $p|F_{p+1}$, then $R_{2k}^{p+1}\equiv -I_{2k}\pmod p$ and
$R_{2k+1}^{p+1}\equiv I_{2k+1}\pmod p$. In our case, $3|F_4=3$.
Thus, the minimal polynomial of $R_n$ modulo 3 equals $x^4+1$ for $n$
even and $x^4-1$ for $n$ odd.
Since $\gcd(3,4)=1$, we deduce that $x^4\pm 1$ has distinct roots.
Thus, $R_n$ is diagonalizable (over a degree 2 extension of ${\bf Z}_3$).

First, $n=4k$. Thus, the minimal polynomial of $R_n$ modulo 3
equals $x^4+1=(x^2-x-1)(x^2+x-1)$, which has the solutions
$\alpha,\beta,-\alpha,-\beta$ (where $\beta=1-\alpha$)
in the quadratic extension of ${\bf Z}_3$,
namely ${\bf Z}_3[x]/(x^2+x-1)$. Let $a,b,c,d$ be the multiplicities of
$\alpha,\beta,-\alpha,-\beta$, respectively.
We want to show that $a=b=c=d=k$.
Since $R_n$ has dimension $n=4k$,
\[
a+b+c+d=4k.
\]
Now, using Lemma \refth{trace}, we get
\[
a\alpha+b\beta-c\alpha-d\beta\equiv F_{4k}\pmod 3.
\]
It is known (see \cite{GKP}) that $F_n|F_{nk}$. Thus, $F_{4k}$ is
divisible by $F_4=3$, so $F_{4k}\equiv 0\pmod 3$. Therefore, $a=c$
and $b=d$, so $a+b=c+d=2k$. The eigenvalues of $R_{4k}$ correspond
to eigenvalues of $R_{4k}^3$ and since, $R_{4k}^3=-R_{4k}^{-1}$ (modulo 3),
we get that $\alpha$ corresponds to $\alpha^3$, which corresponds
to $-\frac{1}{\alpha}=\beta$, and similarly,
$\displaystyle\beta\Longrightarrow
\beta^3\Longrightarrow -\frac{1}{\beta}=\alpha$. We deduce $a=b$
and so, $a=b=c=d=k$.

In a similar manner we can show the other cases.
\end{proof}

With the same method, we can prove
\begin{thm}
The characteristic polynomial of $R_n$ modulo $5$, say $p_n(x)$, is
\begin{eqnarray*}
& p_{4k}(x)=(x-2)^{4k}; & p_{4k+1}(x)=-(x-1)^{4k+1}\\
& p_{4k+2}(x)=(x+2)^{4k+2}; & p_{4k+3}(x)=-(x+1)^{4k+3}
\end{eqnarray*}
\end{thm}

The interesting fact is  that we were able to find,
after quite a bit of work, the eigenvalues (with their multiplicities),
for all the first few primes $p$ we considered.

\section{Characteristic Polynomial, Eigenvalues and Eigenvectors of $R_n$}

In \cite{PS}, the authors proposed a conjecture on
the eigenvalues.
No pattern leaped out from the eigenvectors' set, if one
looks at the first few values.
In this section we prove the conjecture and we
find the eigenvectors of $R_n$.
In \cite{PS}, we showed 
\begin{lem}
The inverse of $R_n$ is the matrix
 \[
 R_n^{-1}=\left((-1)^{n+i+j+1}\binom{n-i}{j-1}\right)_{1\leq i,j\leq n}.
 \]
\end{lem}

We define $K_n$ to be the matrix with $(i,j)$-entry $\delta_{i,n-j+1}$
(the Kronecker symbol).
\begin{rem}
We remark that $K_n$ is the permutation matrix having $1$ on the secondary
diagonal and $0$ elsewhere.
\end{rem}

\begin{thm}
Denote $\phi=\frac{1+\sqrt{5}}{2},\bar\phi=\frac{1-\sqrt{5}}{2}$.
The eigenvalues of $R_n$ are:\\
1. $\displaystyle\{ (-1)^{k+i} \phi^{2i-1},
(-1)^{k+i} \bar\phi^{2i-1}  \}_{i=1,\ldots,k}$, if $n=2k$.\\
2. $\displaystyle\{ (-1)^k \}\cup\{ (-1)^{k+i} \phi^{2i},(-1)^{k+i}
\bar\phi^{2i}  \}_{i=1,\ldots,k}$, if $n=2k+1$.
\end{thm}
\begin{proof}
First, we show that $R_n$ is a permutation
matrix away from $L_n$, namely
\begin{equation}
\label{equiv}
R_n\cdot K_n=L_n\Longleftrightarrow R_n=L_n\cdot K_n.
\end{equation}
It is a trivial matter to prove $K_n^2=I_n$, which will give the
equivalence. It suffices to show the second identity, which follows
easily, since an entry in $L_n\cdot K_n$ is
\[
\begin{split}
\sum_{k=1}^n \binom{i-1}{k-1}\delta_{k,n-j+1}=\binom{i-1}{n-j}.
\end{split}
\]
Now, denote by $A_n$, the matrix obtained by taking
 absolute values of entries of $R_n^{-1}$.
 We show that $R_n$ has the same characteristic
 polynomial (eigenvalues) as $A_n$, namely we prove
 their similarity,
\begin{equation}
\label{sim}
K_n\cdot R_n\cdot K_n=A_n.
\end{equation}
Since $K_n\cdot R_n\cdot K_n=K_n\cdot L_n$ (by \refeq{equiv}), to show \refeq{sim}
it suffices to prove that
$K_n\cdot L_n=
\left ( \binom{n-i}{j-1} \right)_{i,j}$. Therefore,
we need
\[
\sum_{k=1}^n \delta_{i,n-k+1}\binom{k-1}{j-1}=\binom{n-i}{j-1},
\]
which is certainly true.

We use a result of \cite{LP} to show that $A_n$ in turn is similar to
$D_n$, the diagonal matrix whose diagonal entries are the elements
in the eigenvalues set listed in decreasing order according to size
of the absolute value. For instance, for $n=4$,
\[
D_4=
\left(\begin{array}{cccc}
\alpha^3 &0  & 0 & 0   \\
0 & -\alpha & 0 & 0    \\
0 & 0 & -\beta & 0   \\
 0& 0 & 0 & \beta^3
\end{array}\right)
\]
If $A_n$ is similar to $D_n$, the theorem will be proved.
We sketch here the argument of \cite{LP} for
the convenience of the reader.
Define an array $b_{n,m}$ by
\[\begin{split}
b_{n,0}&=1\ \text{for all}\ n\geq 0\\
b_{n,m}&=0\ \text{for all}\ m>n\\
b_{n,m}&=b_{n-1,m-1}\, (-1)^m\, \frac{F_n}{F_m}\ \text{for all}\ m\leq n
\end{split}
\]
Except for the sign,
$\displaystyle
|b_{n,m}|=\frac{F_nF_{n-1}\cdots F_{n-m+1}}{F_m\cdots F_1}$.
Let $C_n=(c_{i,j})_{i,j}$, where
\[
\begin{cases}
c_{i,i+1}=1& if i=1,\ldots, n-1\\
c_{n,j}= -b_{n,n+1-j} & if j=1,\ldots, n\\
c_{i,j}=0& otherwise.
\end{cases}
\]
We observe that, in fact, $C_n$ is the companion matrix of the polynomial
with coefficients $b_{n,n+1-j}$.
Let $X_n$ be the matrix with entries
$\binom{n-i}{j-1} F_{i-2}^{j-1} F_{i-1}^{n-j}$.
%and $M_n$ be the matrix with entries
%$\binom{n-i}{j-1} F_{i-1}^{j-1} F_{i}^{n-j}$.
It turns out that the eigenvector matrix $E_n$ of $A_n$,
with columns vectors listed in decreasing order of
absolute value of the corresponding eigenvalues, normalized so that
the last row is made up of all 1's, satisfies
\[
X_n E_n=V_n,
\]
where $V_n$ is the Vandermonde matrix, which is the eigenvector
matrix of $C_n$ with eigenvectors listed in decreasing order of the absolute
values of the corresponding eigenvalues.
Also,
\[
X_n A_n X_n^{-1}=C_n\ \text{and}\ E_n^{-1} A_n E_n=D_n.
\]
Our theorem follows.
\end{proof}

Easily we deduce
\begin{cor}
The eigenvectors matrix of $R_n$, say $W_n$, with eigenvectors listed in
decreasing order of the absolute
values of the corresponding eigenvalues, is
\[
W_n=K_n E_n=K_n X_n^{-1} V_n.
\]
\end{cor}
\begin{proof}
We showed that
\[
K_n R_n K_n=A_n\ \text{and}\ E_n^{-1} A_n E_n=D_n.
\]
It follows
that $(E_n^{-1} K_n) R_n (K_n E_n)=D_n$, which together with $X_n E_n=V_n$,
proves the corollary.
\end{proof}

\begin{example}
For $n=4$, the eigenvectors matrix is
\[
W_4=
\left(\begin{array}{cccc}
-\alpha^3 &\alpha  & \beta & -\beta^3   \\
\alpha^2 & -\frac{1}{3}\beta & -\frac{1}{3}\alpha & \beta^2    \\
-\alpha & -\frac{1}{3}\alpha^2 & -\frac{1}{3}\beta^2 & -\beta   \\
 1& 1 & 1 & 1
\end{array}\right)
\]
\end{example}

\noindent{\bf Acknowledgements.}
{
The authors would like to thank Richard P. Stanley for providing the references
\cite{Strauss, Waterhouse}.
}

\end{document}